\documentclass[reqno]{amsart}

\usepackage{graphicx}
\usepackage{amssymb}
\usepackage{amstext}
\usepackage{amsthm}
\usepackage{epstopdf}
\newcommand{\R}{{\mathbb R}}

\newtheorem{theorem}{Theorem}
\newtheorem{corollary}[theorem]{Corollary}

\newtheorem{proposition}[theorem]{Proposition}
\newtheorem{conjecture}{Conjecture}
\newtheorem{definition}{Definition}
\newtheorem{example}{Example}

\theoremstyle{remark}
\newtheorem{remark}{Remark}

\numberwithin{equation}{section}
\begin{document}

\title[delay differential equations with negative Schwarzian]{Dichotomy results for  delay
differential equations with negative Schwarzian}
\author[E. Liz and G. R\"{o}st]{}
%\email{eliz@dma.uvigo.es}
%\email{rost@math.u-szeged.hu}
%\subjclass{34K20, 34D45}
\maketitle

\centerline{\scshape  Eduardo Liz\footnote{Corresponding author; E-mail address: eliz@dma.uvigo.es, Fax: (+34)
986 812116.}} \medskip

{\footnotesize
\centerline{ Departamento de Matem\'atica Aplicada II,
E.T.S.I. Telecomunicaci\'on }
\centerline{ Universidad de Vigo, Campus
Marcosende, 36310 Vigo, Spain} }

\medskip \centerline{\scshape  Gergely R\"{o}st\footnote{E-mail address: rost@math.u-szeged.hu}} \medskip

{\footnotesize
\centerline{ Analysis and Stochastics Research Group,
Hungarian Academy of Sciences}
\centerline{ Bolyai Institute, University of
Szeged, H-6720 Szeged, Aradi v\'ertan\'uk tere 1., Hungary} }

\medskip \centerline{\today}

\bigskip

\begin{quote}
{\normalfont\fontsize{8}{10}\selectfont{\bfseries Abstract.} We gain further insight into the use of the Schwarzian derivative to obtain new results for a family of functional differential equations including the famous Wright's equation and the Mackey-Glass type delay differential equations.
We present some dichotomy results, which allow us to get  easily computable bounds of the
global attractor. We also discuss related conjectures, and formulate new open
problems. }
\end{quote}

\section{Introduction}

The property of negative Schwarzian derivative plays a prominent role in the
qualitative analysis of discrete dynamical systems since the famous paper of
Singer \cite{si}.  This fact has recently been used to analyze the behaviour of the
solutions of a family of delay differential equations, in particular to get sufficient conditions for the global stability of the equilibrium. Our aim in this paper is to gain further insight into the use of this property to obtain
some easily computable bounds for the
global attractor of the delay differential equation
\begin{equation}  \label{1}
x^{\prime }(t)=-ax(t)+f(x(t-\tau)),
\end{equation}
where $a\geq 0$, $\tau>0$, and $f$ is a continuous function.
Some of our main conclusions are given in the form of dichotomy results;  related work for a class of cyclic systems may be found in \cite{en}. 
We also discuss some
conjectures concerning the global behaviour of the solutions, suggested by
the analogy with one-dimensional maps, and some numerical experiments.  In the particular case of decreasing maps, we believe that the simple dynamics of the map $f$ given in Theorem \ref{tsinger} below  is inherited by the delay differential equation (\ref{1}). This is stated as Conjecture \ref{co1} in Section 4.

Equation (\ref{1}) has been used as a mathematical model in many different
areas, including models in neurophysiology, metabolic regulation, optics and
agricultural commodity markets. See, e.g., the interesting list in \cite[%
p.~78]{hm} and the introduction in \cite{mp2}.
As noticed in \cite{gtb, gt, lprtt, ltt,mp2}, in some of the most famous models as the Nicholson's blowflies
equation \cite{gbn}, the Mackey-Glass equation \cite{mg}, and Wright's equation \cite{wr}, the involved nonlinearity $f$ is either
decreasing or unimodal and has negative Schwarzian derivative. 

The paper is organized as follows: in Section 2 we review some standard
definitions and properties of one-dimensional maps - in particular those with
negative Schwarzian derivative -, and some basic facts about the delay differential equation (\ref{1}).   In Section 3, we prove the main results for (\ref{1}), assuming that $f$ is decreasing and has negative Schwarzian derivative. Our results are illustrated with examples for the equations of Wright and Mackey-Glass. Section 4 is devoted to discuss some conjectures on the global  behaviour of the solutions of Eq. (\ref{1}). Finally, in Section 5 we formulate an open problem for the case of unimodal feedback.

\section{Preliminaries}
\subsection{One-dimensional maps}
First, we recall some basic concepts about one-dimensional maps, introduce some necessary definitions, and state an easy but important result for our purposes.

By a one-dimensional map we mean a continuous function $f:I\rightarrow I$,
where $I$ is a real interval. For each $n\in {\mathbb{N}}$ we denote
by $f^{n}$ the corresponding power of $f$ under composition, that is,
\[
f^{n}=\underbrace{f\circ \cdots \circ f}_{n}.
\]

To each $x\in I$ we can associate the orbit of $x$, given by the set

\[
\{f^n(x)\, :\, n=0,1,2,\dots\}=\{x,f(x), f^2(x),\dots \}.
\]

A fixed
point $K\in I$ of $f$  is attracting if 
\begin{equation}
\label{attract}
\lim_{n\rightarrow \infty }f^{n}(x)=K,
\end{equation}
 for all $x$ in some neighbourhood of $K$. Assuming that $f$ has a unique fixed point $K$,   we say that $K$ is globally attracting if (\ref{attract}) holds for all $x\in I$.
 For any differentiable function $f$, it
is well known that $K$ is attracting if $|f^{\prime }(K)|<1$.

A set $X\subset I$ is called invariant if $f(X)\subset X$. An invariant
subset $X$ is attracting if there is a neighbourhood $U$ of $X$ such that $%
\bigcap_{n=0}^{\infty}f^n(U)\subset X$.

A set $\{\alpha ,\beta \}$ such that $f(\alpha )=\beta $, $f(\beta )=\alpha $
is called a cycle of period $2$ of $f$ (or a $2$-cycle). Clearly, a $2$%
-cycle is invariant. We say that the $2$-cycle $\{\alpha ,\beta \}$ is
globally attracting if $f^{n}(x)\rightarrow \{\alpha ,\beta \}$ as $%
n\rightarrow \infty $ for all $x\neq K$.

Throughout the paper, we shall consider  maps with negative Schwarzian derivative. We recall
that the Schwarzian derivative of a $C^3$ map $f$ is defined by the relation
\[
(Sf)(x)=\frac{f^{\prime \prime \prime }(x)}{f^{\prime }(x)}-\frac{3}{2}\left(%
\frac{f^{\prime \prime }(x)}{f^{\prime }(x)}\right)^2,
\]
whenever $f^{\prime }(x)\neq 0$.

\begin{definition}
\label{smap}
We say that $f:I\to I$ is an $S$-map if it is three times
differentiable, $(Sf)(x)<0$, and  $f^{\prime }(x)<0$ for all $x\in I$. If $I=\R$, we also assume that 
at least one of the limits $f(-\infty)$ and $f(\infty)$ is finite.
\end{definition}
We emphasize that an $S$-map $f$ always has a compact invariant and attracting interval $[A,B]$. For example, if $f(-\infty)=c\in\R$, then we can choose $B=c$, $A=f(c)$.

\begin{definition}
\label{sumap}
We say that a $C^{3}$ map $f:[a,b]\to [a,b]$ is an $SU$-map if it has a unique critical point $x_0$,  with $f^{\prime }(x)>0 $ for $x<x_0$, 
$f^{\prime }(x)<0$ for $x>x_0$,  
and $(Sf)(x)<0$ for all $x\neq x_{0}$. Moreover, we assume that there is a unique fixed point $K$ of $f$, and $K>x_0$.
\end{definition}

The following dichotomy result for $SU$-maps with negative Schwarzian derivative based on Singer's paper \cite{si}  will be very important in our discussions:
\begin{theorem}
\label{tsinger} Assume that $f:I\to I$ is either an $S$-map or an $SU$-map, and let $K$
be the unique fixed point of $f$ in $I$. Then exactly one of the following
occurs:

\begin{itemize}
\item[\textrm{(a)}] If $|f^{\prime }(K)|\leq 1$, then 
$\lim_{n\to\infty}f^{n}(x)=K\, ,\;\forall\, x\in I.$

\item[\textrm{(b)}] If $|f^{\prime }(K)|>1$, then there exists at most one attracting periodic orbit of $f$. Moreover, if $f$ is an $S$-map, then this orbit exists, has period two, and is globally attracting.
\end{itemize}
\end{theorem}
\proof  A detailed proof of part (a) may be seen in \cite[Proposition 3.3]{ltt}. The existence of at most one attracting periodic orbit follows from Theorem 2.7 in \cite{si}. Finally, assume that $f$ is an $S$-map and $|f'(K)|>1$. Thus, since $(f^2)'(K)=(f'(K))^2>1$,  it follows  from  Singer's results and the existence of an attracting invariant interval $[A,B]$ that the increasing map $f^2$   has exactly three fixed points $\alpha<K<\beta$. Moreover, $f(\alpha)=\beta$, 
$f(\beta)=\alpha$, and
\[
\lim_{n\to\infty}f^{2n}(x)=\left\{
\begin{array}{l}
\alpha\, \text{ if }\; x<K \\
\noalign{\medskip} \beta\, \text{ if } \; x>K.%
\end{array}
\right.
\]
This means that $\{\alpha,\beta\}$ is globally attracting for $f$.
\qed

Thus, the global dynamics of an $S$-map is very simple; taking $%
\lambda=f^{\prime }(K)$ as a parameter, there is a unique bifurcation point $%
\lambda=-1$, which leads to a period-doubling bifurcation. Moreover, the
fixed point is a global attractor for $\lambda\in [-1,0)$, and, for $%
\lambda<-1$, the equilibrium becomes unstable and there is a globally
attracting $2$-cycle. 
For $SU$-maps, this is not true in general.
 In fact, the discrete dynamical system generated by $f$ can
exhibit chaotic behaviour. For example, see the analysis made for
the logistic family $f_{\lambda}(x)=\lambda x(1-x),\,
0\leq \lambda\leq 4,$ in \cite{is}, or for the family of unimodal maps $%
f_{\lambda}(x)=\lambda x e^{-x},\, \lambda>1,$ in \cite{se}.
On the other hand, there are examples of $SU$-maps without attracting cycles; moreover, this fact typically leads to the existence of an ergodic absolutely continuous measure for $f$ (see, e.g., \cite{ce}). 

We notice that  the dynamics  for $S$-maps still  holds for $SU$- maps if the additional assumption $f^2(x_0)\geq x_0$ is satisfied, where $x_0$ is the unique critical point of $f$. Indeed, in this case $f$ is an $S$-map in the attracting interval $[f^2(x_0), f(x_0)]$.

\subsection{Delay differential equations}
Next we recall some basic facts about the delay differential equation (\ref{1}) with $a\geq 0$ and $f: I \to I$, where either $I=\R$ or $I=[0,\infty)$.  We assume that $f$ is continuously 
differentiable. (A good source to read more details is \cite{kr2}.)

Let $C$ denote the Banach space of continuous functions $\phi: [-\tau, 0] \to\R$ with 
the usual sup-norm defined by $\|\phi\| = \max\{ |\phi(t)|\, :\, -\tau\leq t\leq 0\}. $
If $x$ is a solution of (\ref{1}), then the segment $x_t \in C$ is defined by 
$x_t (s) = x(t + s), -\tau \leq s \leq 0.$

Every $\phi \in C$ uniquely determines a solution $x = x^{\phi} : [-\tau, \infty) \to\R$ of (\ref{1}), 
i.e., a continuous function $x: [-\tau, \infty) \to\R$ such that $x$ is differentiable on $(0,\infty),$ 
$x_0 = \phi$, and $x$ satisfies (\ref{1}) for all $t > 0$. 

 The solution map $F:\R^+\times C\to C$ defines a continuous semiflow. 
If $\phi\in C$ and $x^{\phi}$ is bounded on $[-1,\infty)$, then the $\omega$-limit set 
\begin{align*}\omega(\phi) = &\{\psi \in C : \mbox{ there is a sequence $(t_n )$ in $[0, \infty)$ 
with $t_n \to\infty$}\\
&\mbox{ and $F (t_n , \phi) \to\psi$ as $n \to\infty$}\} 
\end{align*}
is nonempty, compact, connected and invariant. According to the Poincar\'e-Bendixson theorem of Mallet-Paret and Sell \cite{mps}, if $f$ is an $S$-map then the $\omega$-limit set of any solution of (\ref{1}) is either the fixed point or a single, nonconstant, periodic solution.

The global attractor of the semiflow $F$ is a nonempty compact set $A \subset C$
which is invariant in the sense that $F (t, A) = A$ for all $t \geq 0$, and which attracts 
bounded sets in the sense that for every bounded set $B \subset C$ and for every open neighbourhood $U$ of $A$ there exists $t \geq 0$ with $F ([t, \infty) \times B) \subset U.$ 
The global attractor is characterized by the following property: 
$$A = \{\phi \in C \,:\, \mbox{there is a bounded solution $x: \R \to\R$
of Eq. (\ref{1}) so that $x_0 = \phi$}\}. $$
If $f$ is an $S$-map in the sense of Definition \ref{smap}, then the global attractor $A$ exists. The same applies to $SU$-maps if $a>0$.
In Section 3, we  get some bounds for the global attractor $A$ of (\ref{1}). To be more precise, saying that an interval $[\alpha,\beta]$ contains the
global attractor we mean that for each $\phi \in A$, we have $\alpha
\leq \phi(s) \leq \beta$ for any $s \in [-1,0].$

\section{Dichotomy results for  monotone feedback}

In this section, we consider Eq. (\ref{1}) assuming that $f$ is an $S$-map. In this case, Theorem \ref{tsinger} may be used to get some dichotomy results. We distinguish the cases $a=0$ and $a>0$. Having in mind the most famous examples, we will refer to the first one as {\it Wright type equation}, and to the second one as 
 {\it Mackey-Glass type equation} .

\subsection{Wright type equations}
First we consider Eq. (\ref{1}) with $a=0$, that is,
\begin{equation}  \label{10}
x^{\prime }(t)=f(x(t-1)),
\end{equation}
where $f:\R\to \R$ is an $S$-map. Without loss of generality, we assume that $f(0)=0$. We also notice that it is not restrictive to assume $\tau=1$, since this may also be achieved by rescaling the time.

 It is very well-known that in this case all solutions of (\ref{10}) with continuous initial condition $x(t)=\varphi(t)$ for all $t\in [t_0-1,t_0]$ are bounded and defined in $[t_0-1,\infty)$.  In this case, we can prove the following dichotomy result:

\begin{proposition}
\label{prw} Assume that $f:\R\to \R$ is an $S$-map and $f(0)=0$. Then:

\begin{itemize}
\item[\textrm{(a)}] If $|f^{\prime }(0)|\leq 3/2$, then $\lim_{t\to%
\infty}x(t)=0$ for every solution $x(t)$ of (\ref{10}).

\item[\textrm{(b)}] If $|f^{\prime }(0)|>3/2$, then
\[
\alpha\leq\liminf_{t\to\infty}x(t)\leq \limsup_{t\to\infty}x(t)\leq\beta,
\]
for every solution $x(t)$ of (\ref{10}), where $\{\alpha,\beta\}$ is the
unique $2$-cycle of $f$.
\end{itemize}
\end{proposition}
\proof
Part (a) is a direct corollary of Theorem 1.3 in \cite{lprtt}. To prove part (b), we recall that every solution of (\ref{10}) is oscillatory if $|f'(0)|>1$ (see, e.g., \cite[Proposition XV.2.2]{detal}). Let $x(t)$ be a solution of (\ref{10}),
 $M = \limsup_{t \to \infty}x(t), \ m = \liminf_{t \to
\infty}x(t)$, $m\leq 0\leq M$. 
There are two sequences of points
$t_j, s_j$ of local positive maxima  and local negative minima, respectively, such
that $x(t_j)= M_j \to M, x(s_j) = m_j \to m$ and $s_j, t_j \to
+\infty$ as $j \to \infty$.
Moreover, for each $t_j$
we can find $\varepsilon_j \to 0+$ such that $m - \varepsilon_j <x(s)$ for all $s\in [t_j-2,t_j-1]$. 
Since $x'(t_j)=0=f(x(t_j-1))$, it follows that $x(t_j-1)=0$. Hence, since $f$ is decreasing,
$$
M_j =x(t_j)= \int_{t_j-1}^{t_j}f(x(s-1))ds <f(m-\varepsilon_j).
$$
As a limit form of this inequality, we get $M \leq
f(m)$.

An analogous argument proves that $m\geq f(M)$, and therefore 
$m\geq f(M)\geq f^2(m)$, $M\leq f(m)\leq f^2(M)$. According to Theorem \ref{tsinger}, this is only possible if $[m,M]\subset [\alpha,\beta]$.
\qed

In the applications, finding the exact values of $\alpha$ and $\beta$ might be difficult. 
Assuming that $f(-\infty)=r$, it is clear that $[\alpha,\beta]\subset [f(r),r]$. This interval may be sharpened to $[f^{2k+1}(r), f^{2k}(r)]$ for any $k=1,2,\dots$.

A famous example of (\ref{10}) with $S$-map $f$ is 
\begin{equation}
\label{wre}
x'(t)=-r\left(e^{x(t-1)}-1\right),
\end{equation}
which is obtained from  the famous Wright equation
\begin{equation}  \label{wright}
y^{\prime }(t)=-r y(t-1) (1+y(t)),
\end{equation}
where $r>0$. This equation was investigated in an outstanding paper
published by E. M. Wright in 1955 \cite{wr}. 
Corresponding to an initial
function $\phi$ such that $\phi(0)>-1$, there is a unique solution $y:[-1, \infty)\to\R$; moreover, $y(t)>-1$ for all $t>0$. Hence, if we only consider this set of initial functions,  the change of variables $x(t)=\ln(1+y(t))$ transforms (\ref%
{wright}) into (\ref{wre}).

It is clear that $f(x)=-r(e^x-1)$ is an $S$-map, $f(-\infty)=r$, and $f'(0)=-r$. 
An immediate consequence of Proposition \ref{prw} is that all solutions of (\ref{wright}) with $y(0)>-1$ converge to $0$ if $r\leq 3/2$, and, if $r>3/2$ and $y(t)$ is a solution of  (\ref{wright}) with $y(0)>-1$, then
\begin{equation}
\label{bound}
-1+\exp\left(-r(e^r-1)\right)\leq  \liminf_{t \to \infty}y(t)\leq  \limsup_{t \to \infty}y(t) \leq e^r-1.
\end{equation}
The first part is the famous {\it Wright's $3/2$- global stability theorem} (see \cite[Chapter 4]{kuang}). Wright also found the upper bound 
in (\ref{bound}). We emphasize that the lower bound is of independent interest, since it shows that the solutions are bounded away from $-1$. This is a result of uniform persistence that is important from the biological point of view (we recall that Wright's equation is equivalent to the delayed logistic equation proposed by Hutchinson via the change of variables $N(t)=1+y(t)$).

However, the upper bound $e^r-1$
 is not very sharp, as we show in Example \ref{ex1} below. 
 In order to improve it, we use a different $S$-map.

 \begin{proposition}
\label{prw2} Assume $r>1$, and define the map
\[
F(y):=-1+\exp\left[\left(ry+1-e^{ry}\right)/y\right].
\]
Let $y(t)$ be a solution of (\ref{wright}) with $y(0)>-1$, and denote $%
m=\liminf_{t\to\infty} y(t)$, $M=\limsup_{t\to\infty} y(t)$.  Then:

\begin{itemize}
\item[\textrm{(a)}] If $r\leq 3/2$, then $\lim_{t\to%
\infty}y(t)=0$.

\item[\textrm{(b)}] If $r>3/2$, then
\[
\alpha_1\leq\liminf_{t\to\infty}y(t)\leq \limsup_{t\to\infty}y(t)\leq\beta_1,
\]
where $\{\alpha_1,\beta_1\}$ is the
unique $2$-cycle of $F$.
\end{itemize}
\end{proposition}
\proof
Part (a) was already proved. 
To prove part (b), we use some relations from the proof of Theorem 3 in \cite{wr}. Indeed, inequality (3.9) in page 70 of this reference (valid for $r>1$) is equivalent to the relation $M\leq F(m)$. The relation $m\geq F(M)$ may be obtained in a completely analogous way.
On the other hand,  the following properties of $F$ are easy to check:
\begin{enumerate}
\item[(i)] $F^{\prime }(y)<0$ for all $y\in{\mathbb{R}}$;
\item[(ii)] $(SF)(y)<0$ for all $y\in{\mathbb{R}}$;
\item[(iii)] $\displaystyle\lim_{y\to\infty}F(y)=-1$.
\end{enumerate}
Thus, $F$ is an $S$-map.
Next, since $|F^{\prime}(0)|=r^2/2>1$ for $r>\sqrt{2}=1.4142$, and we are assuming $r>3/2$, it follows from Theorem \ref{tsinger} that $F$ has a globally attracting $2$-cycle $\{\alpha_1,\beta_1\}$. The relations $m\geq F(M)\geq F^2(m)$, $M\leq F(m)\leq F^2(M)$ imply that $[m,M]\subset [\alpha_1,\beta_1].$ \qed

\begin{figure}[hb]
\centering
\includegraphics[totalheight=2.5in]{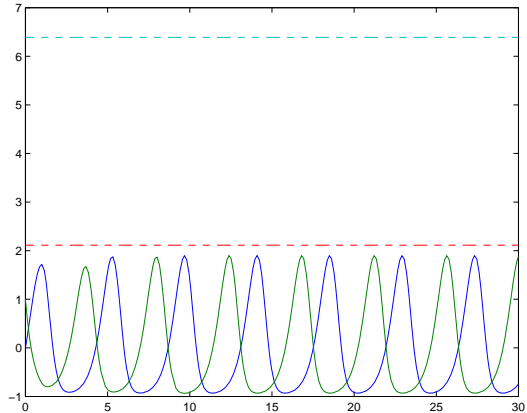}
\caption{Bounds for Wright equation $y^{\prime }(t)=-2 y(t-1) (1+y(t))$.}
\label{wr2}
\end{figure}

As noticed above, the exact values of $\alpha_1$ and $\beta_1$ may be difficult to determine. The following immediate consequence of Proposition \ref{prw2} is useful in the applications:

\begin{corollary}
\label{p4} Assume $r>1$, and let $y(t)$ be a solution of (\ref{wright}) with
$y(0)>-1$. Then,
\begin{equation}
\label{err}
F^2(-1)\leq\liminf_{t\to\infty} y(t)\leq\limsup_{t\to\infty} y(t)\leq F(-1)=-1+\exp\left(-1+r+e^{-r}\right).
\end{equation}
\end{corollary}

The next example shows that the upper bound $F(-1)$ is considerably sharper than the one in (\ref{bound}).

\begin{example}
\label{ex1} \textrm{Consider equation (\ref{wright}) with $r=2$. Using
Proposition \ref{prw}, we get the upper bound $M_1=\beta\approx 6.389$, while Corollary \ref{p4}
provides $M_2=F(-1)=2.112$. See Figure \ref{wr2}, where two solutions of (%
\ref{wright}) are presented, and the horizontal lines indicate $M_1$ and $%
M_2 $ (below).}
\end{example}

\subsection{Monotone Mackey-Glass type equations}
Next we consider the case $a>0$. We can fix $a=1$ to simplify our exposition without loss of generality.
Thus, throughout this subsection we consider the family of delay differential equations
\begin{equation}  \label{mg}
x^{\prime }(t)=-x(t)+f(x(t-\tau)),
\end{equation}
where $f:{\mathbb{R}}^+\to{\mathbb{R}}^+$ is an $S$-map. Although some results may be easily derived from \cite{lr}, we provide here a more systematic approach, and further discussion. 
We mention two
well-known examples of (\ref{mg}) coming from models in hematopoiesis (see, e.g., \cite%
{lptt} and references therein).

\begin{enumerate}
\item The Mackey-Glass equation with decreasing nonlinearity
\begin{equation}  \label{mgd}
x^{\prime }(t)=-x(t)+\frac{p}{1+x^n(t-\tau)},
\end{equation}
with $p>0$, $n>1$.

\item The Lasota-Wazewska equation
\begin{equation}  \label{lw}
x^{\prime}(t)=-x(t)+pe^{-ax(t-\tau)},
\end{equation}
with $p>0$, $a>0$.
\end{enumerate}

We consider only nonnegative solutions of (\ref{mg}). Recall that for each
nonnegative and nonzero function $\phi \in C$,
the unique solution $x^{\phi }(t)$ of (\ref{mg}) such that $%
x^{\phi }=\phi $ on $[-\tau ,0]$ satisfies $x^{\phi }(t)>0\,,\,\forall
\,t>0.$ (See, e.g., \cite[Corollary 12]{gt}.)

Notice that, rescaling the time $t\to t\cdot\tau$, (\ref{mg}) is
equivalent to the singularly perturbed equation 
\begin{equation}  \label{eps}
\varepsilon x^{\prime }(t)=-x(t)+f(x(t-1)),
\end{equation}
where $\varepsilon=1/\tau$. This form was used in \cite{is}, where it was proved
that if $I=[a,b]$ is an invariant and attracting interval for the map $f$,
then $I$ contains the global attractor of (\ref{mg}). 
Combining this result  with Theorem \ref{tsinger}, we have:
\begin{proposition}
\label{c2} Assume that $f$ is an $S$-map and let $K$ be the unique fixed
point of $f$. Then:

\begin{itemize}
\item[\textrm{(a)}] If $|f^{\prime }(K)|\leq 1$, then $\lim_{t\to%
\infty}x(t)=K$ for every solution $x(t)$ of (\ref{mg}).

\item[\textrm{(b)}] If $|f^{\prime }(K)|>1$, then
\[
\alpha\leq\liminf_{t\to\infty}x(t)\leq \limsup_{t\to\infty}x(t)\leq\beta,
\]
for every solution $x(t)$ of (\ref{mg}), where $\{\alpha,\beta\}$ is the
unique $2$-cycle of $f$.
\end{itemize}
\end{proposition}

Proposition  \ref{c2} gives a dichotomy for (\ref{mg}) independent of the delay
$\tau$ which cannot be improved. Indeed, Hadeler and Tomiuk \cite{ht} have proved
that if $|f^{\prime }(K)|>1$ then there exists a $\tau_0>0$ such that (\ref%
{mg}) has a nonconstant periodic solution for all $\tau>\tau_0$. This shows
that part (a) is sharp. On the other hand, part (b) gives the sharpest
invariant and attracting interval containing the global attractor of (\ref%
{mg}) for all values of the delay. This can be proved using Theorems 2.1 and 2.2 in
\cite{mp2} (see \cite[Proposition 5]{lr}).

We emphasize that this interval is quite sharp for sufficiently large
values of the delay, but it is not so accurate for small delays.

\begin{figure}[tbp]
\centering
\includegraphics[totalheight=2.5in]{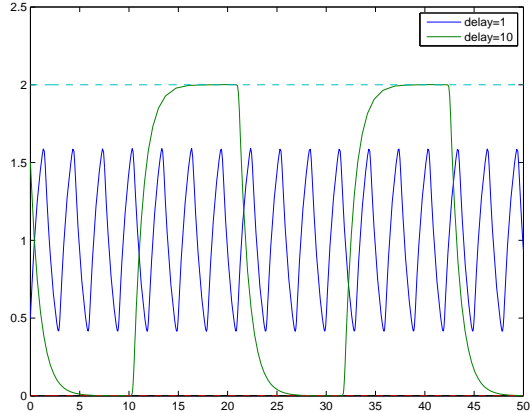}
\caption{Bounds for equation (\protect\ref{mgex2}).}
\label{mgfig}
\end{figure}

\begin{example}
\label{ex2} \textrm{Consider equation (\ref{mgd}) with $p=2, n=20$, that is,
\begin{equation}  \label{mgex2}
x^{\prime }(t)=-x(t)+\frac{2}{1+x^{20}(t-\tau)}\, .
\end{equation}
The equilibrium is $K=1$, and $f^{\prime }(K)=-10$, so that map $f$ has a
globally attracting $2$-cycle $\{\alpha,\beta\}$. One can check that $\alpha$
and $\beta$ are very close to $0$ and $2$ respectively. Thus, interval $%
[0,2] $ contains the global attractor of (\ref{mgex2}) for all values of the
delay $\tau$. We plotted in Figure \ref{mgfig} a solution of (\ref{mgex2})
with $\tau=1$, and other one with $\tau=10$. We can observe that the bound
is much sharper in the latter case. }
\end{example}

Example \ref{ex2} suggests that it would be useful to get bounds for the
global attractor of (\ref{mg}) that depend on the value of $\tau$.

In a very interesting paper \cite{gt}, Gy\H ori and Trofimchuk have improved all
previously known sufficient conditions for the global stability of the
equilibrium of (\ref{mg}), assuming that $f$ has negative Schwarzian
derivative and it is decreasing or unimodal. To achieve this, they have proved an
analogous result of Proposition \ref{c2} for several maps involving the
delay parameter. In particular, we choose the following function due to its
simplicity:
\begin{equation}  \label{defg}
g(x):=\left(1-e^{-\tau}\right) f(x)+e^{-\tau}K.
\end{equation}
The following result is a consequence of the proof of Theorem 14 in \cite{gt}%
:

\begin{proposition}
\label{p6} Assume that $f$ is an $S$-map and let $x(t)$ be a solution of (%
\ref{mg}). If $m=\liminf_{t\to\infty} x(t)$, $M=\limsup_{t\to\infty} x(t)$,
then $[m,M]\subset g([m,M]),$ where $g$ is defined in (\ref{defg}).
\end{proposition}

Notice that $g$ is an $S$-map because $f$ is an $S$-map. On the other hand, $%
|g^{\prime }(K)|=\left(1-e^{-\tau}\right) |f^{\prime }(K)|$. Thus, we have
the following corollary:

\begin{corollary}
\label{c3} Assume that $f$ is an $S$-map and let $K$ be the unique fixed
point of $f$. Then:

\begin{itemize}
\item[\textrm{(a)}] If $\left(1-e^{-\tau}\right)|f^{\prime }(K)|\leq 1$,
then $\lim_{t\to\infty}x(t)=K$ for every solution $x(t)$ of (\ref{mg}).

\item[\textrm{(b)}] If $\left(1-e^{-\tau}\right)|f^{\prime }(K)|>1$, then $%
\alpha_1\leq\liminf_{t\to\infty}x(t)\leq
\limsup_{t\to\infty}x(t)\leq\beta_1, $ for every solution $x(t)$ of (\ref{mg}%
), where $\{\alpha_1,\beta_1\}$ is the unique $2$-cycle of $g$.
\end{itemize}
\end{corollary}

\begin{remark}
As noticed before, sometimes $\alpha_1$ and $\beta_1$ may be difficult to find. Since 
$g^2(0)\leq \alpha_1<\beta_1\leq g(0)$, the interval $[g^2(0),g(0)]$ also
contains the global attractor of (\ref{mg}).
\end{remark}

\begin{figure}[ht]
\centering
\includegraphics[totalheight=2.5in]{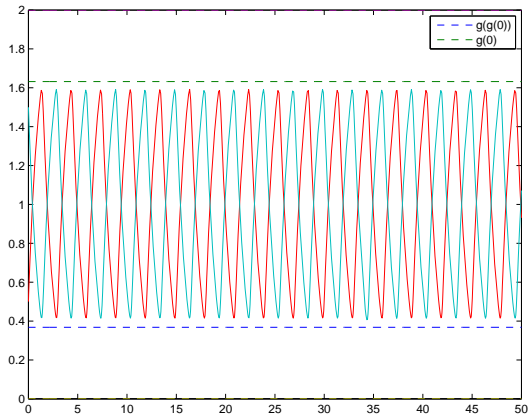}
\caption{Sharper bounds for equation (\protect\ref{mgex2}).}
\label{mgfig22}
\end{figure}

\begin{example}
\label{ex3} \textrm{Consider Eq. (\ref{mgex2}) with $\tau=1$. Then, $%
g(x)=\left(1-e^{-1}\right) f(x)+e^{-1}$, and the interval $[0,2]$ is
sharpened to $[g^2(0),g(0)]=[0.3679,1.6321]$. In Figure \ref{mgfig22} we
depict two different solutions of (\ref{mgex2}), showing that the new bounds
of the global attractor are much more precise than those obtained from the
map $f$. }
\end{example}

\begin{remark}
\label{rgas} We notice that the condition for the global stability of the
equilibrium stated in part (a) of Corollary \ref{c3} has been improved in \cite%
{gt} and later in \cite{lptt, ltt}. However, the expressions of the involved
maps are much more complicated, and this makes difficult to get new easily
computable bounds for the global attractor of (\ref{mg}). For example,
Theorem 14 in \cite{gt} establishes an analogous to Proposition \ref{p6} for
the map $h(x)=F^{-1}\left((1-e^{-\tau})f(x)\right)$, where $%
F(x)=x-e^{-\tau}f^{-1}(x)$. As a consequence, the interval $[h^2(0), h(0)]$
contains the global attractor of (\ref{mg}). In general, it is
difficult to obtain an analytic expression for the inverse functions of $f$
and $F$. However, in the particular case of Example \ref{ex3}, we can easily get $F(x)=x-e^{-1}(2/x-1)^{1/20}$, and we can solve numerically equation $F(h(0))=2(1-e^{-1})$ to improve the bound $g(0)=1.6321$ to $h(0)=1.6071$.
\end{remark}

\begin{remark}
There are particular examples given by Walther \cite{wa2} and Krisztin \cite{kr3}, showing that in Definition \ref{smap} the
conditions $(Sf)(x)<0$ and $f'(x)<0$ are both crucial. In the absence of any of them, these examples show
that not only the dichotomy results of Theorem \ref{tsinger} for the maps are not valid,
but also Propositions \ref{prw} and \ref{c2}: even in case $(a)$ the corresponding delay differential equations
have periodic solutions. 

\end{remark}

\section{On the global dynamics of Equation (\protect\ref{1})}

As noticed in Section 2, if $f$ is an $S$-map then Eq. (\ref{mg}) has a
compact global attractor (see, e.g., \cite{wa}). According to Proposition \ref%
{c2}, the global attractor is $\{K\}$ if $|f^{\prime }(K)|\leq 1.$
On the other hand, if $|f^{\prime }(K)|>1$ then there exists a
positive $\tau _{0}$ such that the equilibrium of (\ref{mg}) is
asymptotically stable if $\tau <\tau _{0}$, and unstable if $\tau>
\tau _{0}$. The value of the bifurcation point $\tau _{0}$ may be
given in terms of $|f^{\prime }(K)|$ (see, e.g., \cite[p.~
135]{hl}).

From Corollary \ref{c3}, it follows that the equilibrium of
(\ref{mg}) is globally attracting for all sufficiently small values
of the delay parameter
$\tau$. As mentioned in Remark \ref{rgas}, this result was improved in \cite%
{ltt}, where it was found a global stability condition very close to the
condition of local stability, that is, for each value of $|f^{\prime }(K)|>1$%
, there is a $\tau_1$ close to $\tau_0$ such that $K$ is globally attracting
for all $\tau<\tau_1$. For this reason, it was conjectured that the local
asymptotic stability of the equilibrium implies its global asymptotic
stability (see more discussions in \cite{ltt2}). 

An analogous conjecture for Wright's equation  (\ref{wre}) 
is very well-known; it states that all solutions of (\ref{wre}) approach zero as $t\to\infty$ if $r\leq \pi/2$. This is usually referred to as the {\it Wright conjecture}, and it is motivated by the paper of Wright in 1955  \cite{wr}. A related conjecture suggested by G. S. Jones  in 1962 \cite{jo} stands for the complementary values of $r$:  \emph{``for every $r> \pi/2$, there exists a unique (up to translations) stable periodic solution of (\ref{wre})"}.
 In this direction, X. Xie \cite{xiet, xie} has proved that the Wright equation has a unique and asymptotically stable slowly oscillating periodic solution when $r\geq 5.67$ (we recall that a periodic solution  is slowly oscillating if the distance between two consecutive zeros is bigger than the delay; otherwise, it is called rapidly oscillating).
This conjecture has recently been revisited by   J.-P. Lessard in his Ph. D. thesis \cite{le}.

Furthermore, both conjectures have been unified and generalized in the recent survey of Tibor Krisztin \cite{kr2}. This generalized conjecture says that the global attractor of (\ref{wre}) is $\{0\}$ if $r\leq \pi/2$, and for $r>\pi/2$ it is formed by zero, a finite number of periodic orbits (the number depends only on $r$), and heteroclinic connections between $0$ and the periodic orbits, and between certain periodic orbits. We call this the {\it Krisztin-Wright conjecture}. It would be interesting to study the same conjecture for Equation (\ref{1})
when $a\geq 0$ and $f$ is an $S$-map. We remark that in the case of $SU$-maps, under some 
additional conditions, all solutions of (\ref{1}) are attracted by the domain where $f'$ is
negative (see \cite{rw} for details), hence analogous conjectures can be 
naturally formulated for a class of delay differential equations with $SU$-maps too.
Motivated by the above discussion and numerical experiments, we guess that the
dichotomy result given for $S$-maps in Theorem \ref{tsinger} may be extended
to the delay equation (\ref{1}). This is stated in the following conjecture:

\begin{conjecture}
\label{co1} Assume that $a\geq 0$ and $f$ is an $S$-map. Let $K$ be the unique fixed
point of $f$. Then:

\begin{itemize}
\item[\textrm{(a)}] The equilibrium $K$ is a global attractor of (\ref{1})
whenever it is asymptotically stable.

\item[\textrm{(b)}] If $K$ is unstable, then there
exists a unique  stable periodic orbit of (\ref{1}), which attracts an open
and dense subset of the phase space.
\end{itemize}
\end{conjecture}

We notice that  much is known about the global dynamics of (\ref{1})
when $f$ is decreasing and bounded from below or from above. A good source is the mentioned paper of Krisztin \cite{kr2},
where the most relevant results and references are included. In particular, Mallet-Paret and Sell \cite{mps} have proved that rapidly oscillating periodic solutions are unstable. On the other hand, combining the main results of Walther \cite{wa}, and Mallet-Paret and Walther  \cite{mpw}, it follows that  if $f$ is $C^1$ and  there is a unique slowly oscillating periodic solution, then  it attracts a dense and open subset  of the phase space $C$. Thus we can rewrite the second part of Conjecture \ref{co1} in the following equivalent form:
\begin{itemize}
\item[\textrm{(b')}] If $K$ is unstable, then Eq. (\ref{1}) has a unique slowly oscillating periodic solution (up to translations in phase).
\end{itemize}

Generalizing the Krisztin-Wright conjecture, we also guess that the global attractor of (\ref{1}), when $f$ is an $S$-map and $K$ is unstable, is formed by $K$, a finite number $N$ of periodic  orbits, and heteroclinic connections between $K$ and the periodic orbits, and between certain periodic orbits.
We notice that $N$ is the number of roots $\lambda$ of the characteristic equation $a+\lambda-f'(K)e^{-\lambda}=0$ associated  to (\ref{1}) with $\mbox{Re}(\lambda)>0$ and $\mbox{Im}(\lambda)>0$ (see \cite{kr2} for details).

Furthermore, if Conjecture \ref{co1} is true, the results of Mallet-Paret and
Nussbaum in \cite{mp1, mp2} ensure that, for $a>0$,  the profiles of the stable periodic
solutions of (\ref{mg}) approach as $\tau\to\infty$ a ``square wave" $p(t)$
defined by the unique $2$-cycle $\{\alpha,\beta\}$ of $f$, that is, $%
p(t)=\alpha$ for $t\in [0,1)$, $p(t)=\beta$ for $t\in [1,2)$, $%
p(t)=p(t+2)$ for all $t$. We emphasize that the limit form of (\ref{mg}) as $\tau\to\infty$ coincides with (\ref{eps}) with $\varepsilon=0$, that is, the difference equation with continuous argument
$$
x(t)=f(x(t-1)).
$$
Thus, the square wave $p(t)$ is defined by the globally attracting $2$-cycle of $f$ given in Theorem \ref{tsinger}.

As noticed by Krisztin (see \cite{kr, kr2} and references therein), Conjecture \ref{co1} is true for a class of odd functions with a convexity property. Two examples are
\[
f_1(x)=-a \tanh(bx)\; ,\quad f_2(x)=-a \tan^{-1}(bx), \quad a>0,\, b>0.
\]
It is not difficult to prove that $f_1$ and $f_2$ are $S$-maps; indeed,
\[
(Sf_1)(x)=-2b^2<0\; ,\quad (Sf_2)(x)=\frac{-2b^2}{(1+b^2x^2)^2}<0,
\quad\forall x\in{\mathbb{R}}.
\]

However, in the general case even the first part of the conjecture seems to
be very difficult to prove (or disprove). We recall that the Wright
conjecture is still an open problem.

\section{Unimodal feedback}

Similar results to that obtained for (\ref{mg}) in the case of $S$-maps may be extended to $SU$-maps; see, e.~g., \cite{lr, rw}. Moreover, the first 
 statement of Conjecture \ref{co1} was also suggested to be true in  this case;
 for the particular case of the Nicholson's blowflies delay differential equation, this problem  was proposed by Smith in \cite{hs}; see \cite{ltt, ltt2} for more discussions.
However, since even chaotic behaviour is possible (see, e.g., \cite{lw}), there is no hope to expect that part (b) hold for $SU$-maps.

We would like to finish suggesting a different question. Taking into account that $SU$-maps have at most one stable periodic orbit, it would be interesting to investigate whether or not a similar property holds for Eq. (\ref{mg}) with $SU$-nonlinearity. If Conjecture \ref{co1} is true, the answer is positive for $S$-maps. Moreover, in this case there would be a unique stable periodic orbit for all values of the delay, which corresponds  either to an equilibrium point or to a slowly oscillating periodic solution.

As mentioned before, Eq. (\ref{mg}) with monotone negative feedback does not have rapidly oscillating stable periodic solutions. 
However, in the general case of (\ref{mg}) with negative feedback, this result is no longer true. Examples of stable rapidly oscillating periodic solutions may be found in the paper of Ivanov and Losson \cite{il}. See also the recent paper of Stoffer \cite{st} and references therein. These examples involve functions $f$ which are constant in some intervals, so they are not $SU$-maps. However, S. Richard Taylor in his Ph. D. thesis \cite{ta} gives a numerical example of the coexistence of two stable periodic solutions in  the Mackey-Glass equation
 \begin{equation*}
x'(t)=-6.15385\, x(t)+\frac{73.8462\, x(t-1)}{1+x^{10}(t-1)}.
\end{equation*}

 This shows that multistability is possible in  (\ref{mg}) with $SU$-maps. We notice that, in this example, one of the periodic solutions is rapidly oscillating.
Thus, we propose to study the following question:

\noindent {\bf Open problem.} Investigate whether or not  coexistence of two slowly oscillating periodic solutions  for (\ref{mg}) with $SU$-maps is possible.

\section*{Acknowledgements}
E. Liz was partially supported by MEC (Spain)  and FEDER,  grant
MTM2007-60679. G. R\"{o}st was partially supported by the
Hungarian Foundation for Scientific Research, grant T 049516. The authors thank Professors
Sergei Trofimchuk and Tibor Krisztin for their valuable comments and suggestions.

\end{document}